**Multiple Forms of Knowing in Mathematics: A Scoping Literature Study**

Hongzhang Xu[3] and Rowena Ball[4]

Mathematical Sciences Institute, Building 145, Science Road, Australian National University, Canberra ACT 2601 Australia

**Emails**

Hongzhang Xu: hongzhang.xu@anu.edu.au

Rowena Ball: rowena.ball@anu.edu.au

**ORCIDs**
Hongzhang Xu: orcid.org/0000-0001-8904-2976
Rowena Ball: orcid.org/0000-0002-3551-3012

**Abstract**

We present a scoping review of published literature on ethnomathematics and Indigenous mathematics as a step towards a goal to decolonize the prevailing Eurocentric view of the provenance of mathematics. Mathematical practices were identified globally from 169 included studies. We map three development stages of ethnomathematical research from 1984 to 2023 and identify 20 categories of Indigenous and traditional cultural activities that evidence mathematical design and expression. We address two challenges of investigating non-Western based mathematics: where to look for mathematical knowledge, and how to decode it from cultural practices. These two hurdles are overcome by cluster analysis of the keywords of included studies. Existing research falls into two categories: I. identification of mathematical concepts used in Indigenous societies, and II. systematizing identified mathematical concepts. Both approaches are essential for research on Indigenous mathematics to flourish, in order to empower Indigenous knowledge holders and deconstruct restrictive colonial boundaries of mathematical knowledge and education.

---

[3] Corresponding author, email: hongzhang.xu@anu.edu.au, phone: +61 2 6125 0727.
[4] Email: rowena.ball@anu.edu.au

**1. Introduction**

Most mainstream mathematicians today believe that pure mathematics is the abstract science of number, pattern, and relationships (Yadav, 2019). Real-world contexts are considered unnecessary and distracting, and mathematics should be presented in a certain form, which is universally valid. However, this definition ignores diverse cultural faces of mathematics and often belies cultural and demographic groups who perform below average in standardized mathematical tests. These omissions and effects have led Bishop (1990) to describe mainstream mathematics as the 'secret weapon of cultural imperialism' (p.51). Mainstream mathematics (here taken as synonymous with 'Western', 'school' or 'British-European mathematics') was taught as a privileged triumph of European civilization (Sturik, 1954), and that view prevails to this day. The contributions of other cultures are played down or denied[1], including those from Egyptian, Sumerian, Indian, Arab, Chinese and other unnamed Indigenous and traditional societies (Closs, 1986).

Despite growing interest in mathematics in Indigenous, non-European, recreational, oral, professional and trades and other non-mainstream contexts, these diverse expressions and applications continue to be marginalized or spurned (Kline, 1982; Skovsmose, 2022; Vandendriessche & Pinxten, 2023). An overarching descriptor for this broad and evolving field is lacking, and the definitions of its components, such as ethnomathematics or Indigenous mathematics, remain contested (Gerdes, 1994; Pais, 2013). What is ethnomathematics, and what should be included? What is Indigenous mathematics? Amidst these debates and questions, there is at least one consensus: a transformative movement has begun towards more diverse and inclusive mathematics.

---

[1] For example, Charles Whish published strong evidence in 1825 that the fundamentals of the calculus, usually credited exclusively to Newton, Liebnitz, and Gregory, were developed by the Kerala School more than two hundred years earlier, notably by Madhava of Sangamagrama and his disciples, and percolated to Europe (Almeida, D. F., & Joseph, G. G. (2004). Eurocentrism in the History of Mathematics: the Case of the Kerala School. *Race & Class*, *45*(4), 45-59. https://doi.org/10.1177/0306396804043866

This transformation will be crucial to the reformation of mathematical education. In some early outcomes of this global trend, decolonization and recontextualization of the current Western curriculum (Garcia-Olp et al., 2022; Paraide et al., 2023) demonstrably have improved students' performance (Fitiradhy et al., 2023; Shahjahan et al., 2021; Supriadi et al., 2023). However, most education-related studies relate primarily to the redesign of classroom activities, course design or other aspects of pedagogy, rather than questioning the actual mathematics content (Rigney et al., 2020). Worldwide, the current mathematical curriculum largely reflects Eurocentric knowledge and ways of thinking (Canevez et al., 2022). One of the perverse impacts of this dominance has been to instil the misconception that Indigenous peoples are somehow intrinsically inept at mathematics. National surveys and international testing programs such as Program for International Student Assessment (PISA), keep the public informed of 'under-achievement' of Indigenous students in mathematics (Klenowski, 2009; Rigney et al., 2020). Such findings both are informed by and help to entrench racist stereotypes and falsehoods (Nhemachena et al., 2020). Although the process of decolonization of science, technology and engineering is (patchily) underway, for mathematics there is an inertial lag (Almeida & Joseph, 2004; Bishop, 1990).

We do need to interrogate the persistent dogma that 'real' mathematics should be context-free and presented in prescribed form, but we also may view mathematics as the low hanging fruit in a wider range of decolonization efforts, precisely because – like art or music – it has developed and is expressed in all human societies. Although mathematics of specific traditional and Indigenous societies has been investigated, cross-cultural comparative studies are limited in scope. Identification of mathematics that is not European-based remains problematic for researchers steeped in the all-pervasive Eurocentric epistemology and curriculum, which limits understanding of Indigenous contributions (Babbitt et al., 2012). To address these gaps, we conducted a scoping study to: (1) identify studies on mathematics in traditional and Indigenous cultures; (2) understand the antecedents and development of studies in this field; (3) investigate Indigenous usages of mathematical concepts and expressions; (4) help find common ground between mainstream and Indigenous mathematics; (5) serve as a resource for educators; and (6) captivate students by showcasing the rich diversity of mathematical expressions across cultures.

## 2. Methods

Scoping studies originally were developed by medical and health researchers to help inform evidence-based methods and processes in the highly variable and complex settings of healthcare (Arksey & O'Malley, 2005; Cooper et al., 2021). Differing in purposes and objectives from narrative reviews and systematic reviews, which aim to collect and curate a certain literature, a scoping review aims to identify and map the extent and nature of the available research evidence on a question, identify emerging discussions of a new issue (Levac et al., 2010), and reinterpret the literature analytically to inform decision-making and future studies (Møgelvang & Nyléhn, 2023; Munn et al., 2018).

To date no scoping study has been carried out on culturally diverse mathematical practices of traditional and Indigenous societies. This scoping study identifies and maps different forms of knowing in mathematics. We outline briefly the four steps of our procedure, based on the framework developed by Arksey and O'Malley (2005), and elaborate these steps and the results in section 3.

**Step 1: Identifying the research questions**

Our research questions, informed by our desired goals, are 1) What does the landscape of research on mathematics of traditional and Indigenous communities look like? 2) What types of cultural activities are enabled or purposed by mathematics?

**Step 2: Identifying relevant studies**

A search string was made using Boolean operators (including OR, AND, NOT, quotation marks, wildcards and brackets): ("indigenous math*" OR "ethnomathematics"). The search string was applied to Scopus[2] on 28th August 2023.

**Step 3: Study selection and charting the data**

In total, 549 records were identified on Scopus then downloaded and screened according to the flow diagram in Figure S1 of Supplement 1. Articles that primarily are concerned with applying ethnomathematics as an educational approach were excluded, as such studies aim to enrich classroom teaching and improve students' performance rather than focus on mathematical knowledges of traditional and Indigenous cultures.

**Step 4: Collating, summarizing and reporting the results**

After screening, 169 English-language articles providing empirical or specific illustration of Indigenous mathematical practices were included. The screening outcomes are presented in Supplements 1, 2 & 3. The articles then were categorized based on the narrative approach proposed by Arksey and O'Malley (2005). This is an interpretative, rather than aggregative, process where 'narrative' refers to synthesizing the studies in a way that tells a story about the overall findings. It enables us to present both chronological (Section 3.1) and topical (Section 3.2) accounts of research on different forms of knowledge in mathematics.

**3. Results**

---

[2] https://www.scopus.com/search/form.uri?display=basic#basic

**3.1 A brief chronicle**

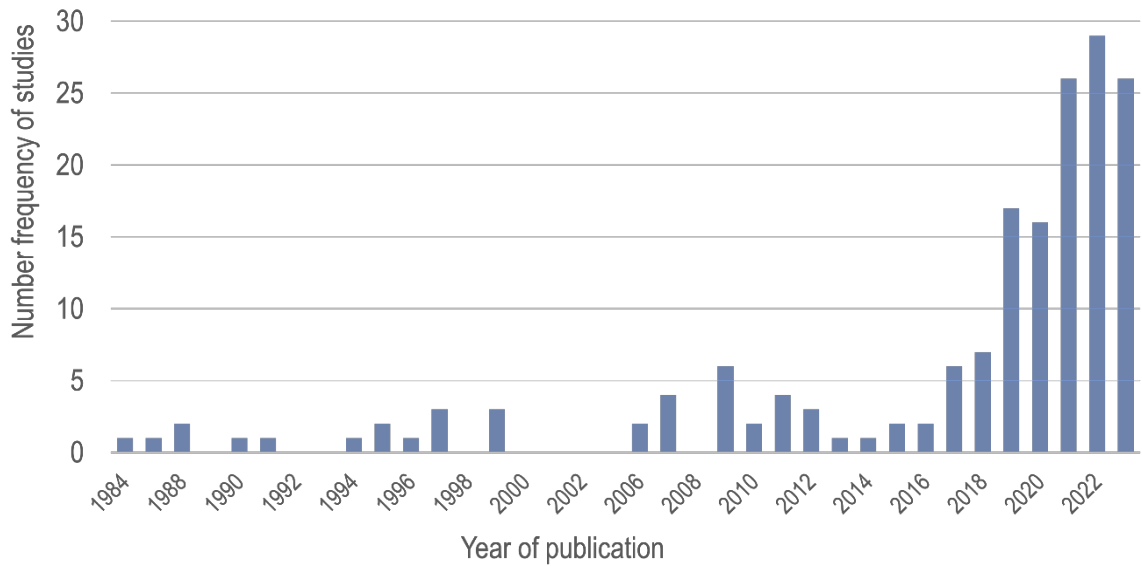

Figure 1. Histogram showing the number of studies of mathematics of traditional and Indigenous communities in each year from 1984 to 2023, for a total of 169 (after screening).

Refer to the histogram of Figure 1. Although research on mathematics of traditional and Indigenous communities can be dated to the 1890s (Mathews, 1900; Spencer & Gillen, 1898), awareness of cultural aspects of mathematics and mathematical education emerged between the late 1970s and the early 1980s (Gerdes, 1994). The first wave of attention pulls together research from two directions: anthropology and education. Mathematical educators found that the history of mathematics presented at junior and high schools tends to downplay or even ignore the cultural context (D'Ambrosio, 1985). The earliest article identified in our scoping study (Lancy, 1981) reported on culturally-aligned disparities in types of mathematical knowledge between seafaring, hunting, and farming Indigenous societies in Papua New Guinea. Ponam and Mandok people performed better than other people on cognitive test because they live on and by the sea, having strong spatial reasoning abilities and good understanding of cycles and patterns (Lancy, 1981).

The recognition of culture-specific expressions of mathematics also inspired new interest in the mathematics developed and used by distinct cultural groups, including occupational labour groups and professional classes, as well as that of Indigenous and traditional peoples. This entire distributed field was called ethnomathematics by D'Ambrosio in the late 1960s (Owens, 2015) and differentiated from mainstream, school or academic mathematics, which was defined as 'mathematics which is taught and learned in schools'. D'Ambrosio (1985) also proposed that some applied mathematics, such that used by practical engineers and applied scientists and technicians, fits better into the category of ethnomathematics mathematics. D'Ambrosio (1995) exposed the colonial origins and values of current schooling and proposed that multicultural and ethnomathematics education is a redeeming opportunity to address inequities in mathematics education and revive traditional cultures.

Anthropologists, represented by Marcia Ascher and Robert Ascher, examined the ways in which non-Western cultures and traditional societies use mathematics, by referring to number, logic and spatial configurations (Ascher & Ascher, 1969, 1986). An alternative definition proposed by Ascher and Ascher (1986) that describes ethnomathematics as studies of mathematical ideas of nonliterate people is also broadly accepted. Their studies of the *quipu*, a medium of ingeniously knotted cords used until the 17th century for record keeping and communication throughout the Inca Empire, are particularly detailed (Ascher & Ascher, 1969, 1972). The use of the *quipu* has been traced back to the Wari Empire (500-1000 AD) and perhaps even to the Caral-Supe civilization (3,000-1,800 BC), both in what is now Peru (UNESCO World Heritage Convention, 2009). Marcia Ascher also identified historical knowledge of graph theory among pre-colonial Indigenous peoples of New Ireland (Niu Ailan), Vanuatu, Angola, Zaire and Zambia, by studying the sand drawing, carving and weaving of complex figures in which a single continuous line is used or implied (Marcia Ascher, 1988; M. Ascher, 1988). Her work belies the predominant belief that graph theory originated from the Königsberg bridge problem studied by Euler (Sachs et al., 1988).

However, opposition to ethnomathematics arose shortly after it gained wider recognition. Chevallard (1990) argued that ethnomathematics was 'enculturation' of mathematics, which in ancient societies such as Babylon and Egypt was 'protomathematics', and that the real mathematics has Greek origins and was passed directly down to form the basis of the present-day Western mathematics corpus. Thomas (1996) dismissed the role of culture in mathematics and emphasized that real mathematics must be decontextualized. Vithal and Skovsmose (1997) and Knijnik (2009) argued that advocates of ethnomathematics did not adequately consider the diverse cultural backgrounds and experiences of students within a school cohort. Non-academic mathematics is also denied recognition as part of the discipline (Chevallard, 1990; Rowlands & Carson, 2002). Some researchers expressed concerns that ethnomathematics could 'pull down' the level of researching and teaching of 'real mathematics'. Rowlands and Carson (2002) decried that mathematics would become a form of 'work experience' (p.84) and declared that 'ethnomathematics runs the risk of attempting to equalise everything down to the poor' (p.98).

Proponents of ethnomathematics have continued to extend their research to diverse cultures and fields. Research articles on exploring mathematics of traditional and Indigenous communities have been increasing in number since 2006 (Figure 1), with ethno-computing and ethno-modelling emerging as prominent areas within the evolving scope of ethnomathematics. Eglash (Eglash et al., 2006; Lotfalian, 2001) proposed the name ethno-computing in 1999 to describe his computer simulations of African Indigenous fractal designs, which he developed into culturally situated design tools (CSDTs), a web-based software, to teach students the mathematical principles in cultural arts such as Native American beadwork, African American cornrow hairstyles and urban graffiti. Rosa and Orey (2010) and Orey and Rosa (2012) proposed that, since mathematics practised by different cultural groups is diverse and distinct, ethno-modelling could be used to help systematize ethnomathematical knowledge and facilitate the understanding of disparate mathematics systems. Uba Umbara et al. (2021) applied an ethno-modelling approach to expose principles such as modulo and congruence applied by the Cigugur Indigenous people in making life choices. Orey and Rosa (2021) set up ethno-models of estimating the volume of tree trunks to elaborate the mathematical practices of local people in Cascavel, Brazil.

Within the body of research on the mathematics of traditional and Indigenous communities, there is a sustained thread on identification of basic geometric motifs and structures, such as cylinders, triangles, and polygons. These studies are influenced by the concept of 'everyday mathematics' from D'Ambrosio (1985) and Lave et al. (1988). We select a few standout exemplars: Stathopoulou (2007) described shapes, colors and patterns of the *xysta* of Pyrgi (on the Greek island of Chios), a type of designed and crafted graffiti on the façades of buildings. Rubio (2016) observed and documented the geometrical patterns made and used by the Kabihug people in sitio Calibigaho, Philippines, within a suite of mathematical activities that included counting, ciphering, measuring, classifying, ordering, inferring and modelling. Utami et al. (2021) identified trapezoid, truncated sphere, tube and other mathematical structures in the designs of Indigenous houses and artifacts in Indonesia. Fauzi et al. (2022) elaborate the space, size and structure in the Sasak traditional architecture and emphasized that mathematics can be found in many activities of daily life.

While these studies offer foundational 'raw materials' for subsequent ethnomathematical research, they risk provoking greater skepticism. This stems from the fact that basic mathematical constructs, such as geometric shapes, can be discerned outside of human societies as well and can be found in nature. Consequently, such a limited focus could inadvertently lead to further de-humanisation of Indigenous and traditional communities, as it may overlook the cultural and complicated ways in which these communities engage with mathematical concepts. Mainstream mathematicians may not recognize this as evidence of Indigenous knowledge of mathematics, but rather as proto-mathematics (Chevallard, 1990).

**3.2 Multiple forms of knowing in mathematics**

Bishop (1991) identifies six 'universal' mathematical activities: counting, locating, measuring, designing, playing and explaining. However, this categorization has not been widely used, perhaps, in part, because ethnomathematical scholarship tends to focus on mathematical ideas within one or a few local cultural groups that are accessible to the scholar (Mukhopadhyay & Roth, 2012; Sizer, 2000). An expanded framework that includes Bishop's categories is developing, strengthened by recognizing that mathematics is an evolving discipline associated with complex cross- and intra-cultural histories.

As an aid to delimiting the preoccupations and restrictions of ethnomathematics research, the 100 most frequent words used by the authors of the 169 included studies are displayed as a wordcloud in Figure S2 of Supplement 1. The word 'culture' (277 occurrences) is the second most frequent and 'tradition' (108 occurrences) is the fifth most frequent word. The word 'Indigenous' has not been broadly used and appeared only 35 times.

Orey and Rosa (2012) state that 'mathematical thinking is influenced by the vast diversity of human characteristics that include: language, religion, worldview, and economical-social political activities' (p.183). All peoples reason about number, logic, and spatial configuration, and organize their ideas into systems. While the European mathematical framework currently is dominant, a deeper dive into history and a broader perspective on knowledge creation and transmission are necessary to gain understanding of mathematical knowledge systems of Indigenous and non-European cultures. (Rosa & Orey, 2010).

As a practical first step toward this objective, we catalogue and classify cultural activities that have obvious direct connections to mathematics. The graphic in Figure 2 visualizes 20 types of cultural activities for which mathematics is required directly. Here we examine many of these categories in more detail.

*Textiles, weaving and carving*

From Figure 2 we see evidence that patterning in textiles and carving (20.5%) is a widespread and well-studied cultural practice that relies on mathematical ways of thinking (Supplements 2 & 3). Weaving crafts and technologies produce textiles with innate patterns that have universal, abstract geometries but culturally specific meanings. Gerdes (1992) and Barton (1996) present a variety of symmetrical patterns of triple weaving that are used in disparate societies, such as in Mozambique, Brazil, China, India, Laos and Japan.

Ornamental arts in many different media have been used to express patterns that have cultural significance. Bérczi (2012) presented the plane symmetry patterns in ancient Eurasian ornamental arts. Rosa and Orey (2008) studied the sacred diamond mat forms rendered in stone, jewelry and cloth and decoded the relations between numbers placed on these mats. More recently, geometric and symmetry transformations used in Indonesian batik patterns were described, including reflection, rotation and translation (Nugraha, 2019; Prahmana & D'Ambrosio, 2020). Faiziyah et al. (2021) expressed reflection-patterning algebraically. However, algebraic methods began to be used in the first half of the 20$^{th}$ century, by 'mathematical weavers'. Ada Dietz (1882-1970) generated weaving patterns based on the polynomial square (O'Brien, 2021)

$$(a + b + c + d + e + f)^2,$$

where the letters represent threads different in colour or texture and each term in the expansion, such as *aa*, *ab* and *ef*, is a certain intertwining of two threads. Building on Ada Dietz's work, Griswold (2002) generated complex, algorithmically-driven designs by incorporating derivatives and integrals of polynomials into weaving patterns that were rendered computationally. Further elaboration of and insight into these algebraic weaving methods is given in Table S1 and Figure S3 of Supplement 1.

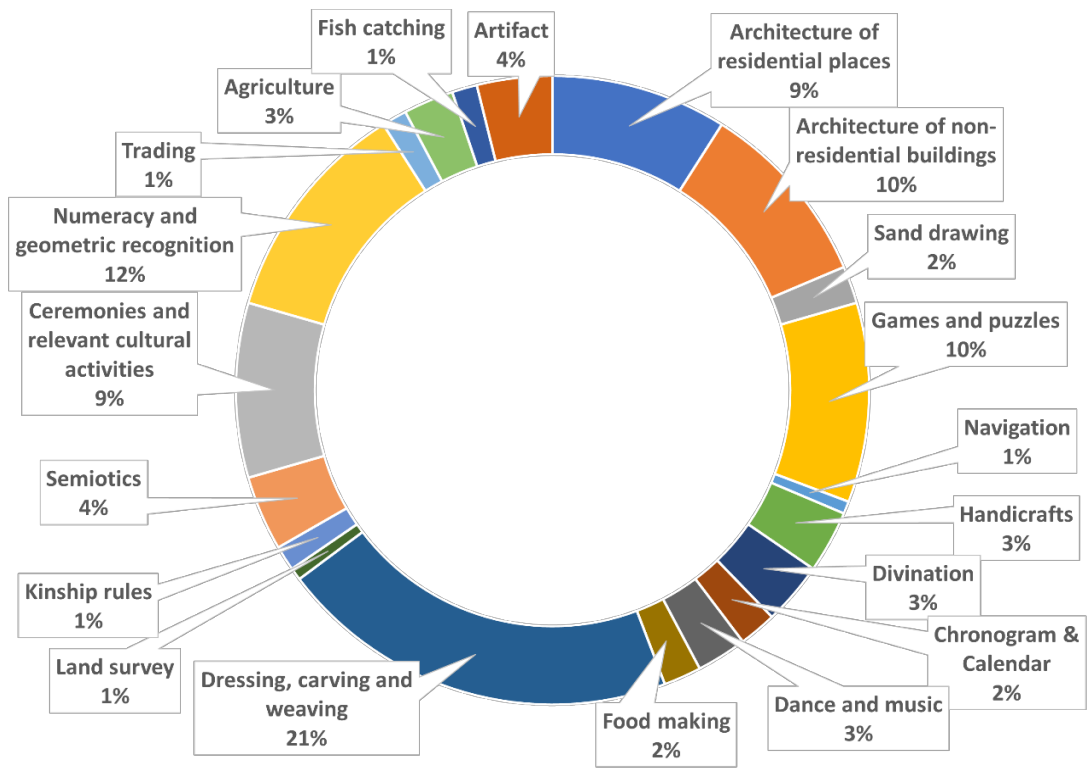

Figure 2. The 20 categories of cultural activities identified in the included studies that involve the explicit use of mathematics. (The total of 102% is due to rounding error.)

*Numeracy and geometric recognition*

Fourteen of the articles (Supplements 2 & 3) describe impressive capabilities in arithmetic, measurement, and estimation in diverse Indigenous communities. Yusharina et al. (2010) describe the volume and weight measurement system of the Indigenous Minangkabau people in Malaysia. Rubio (2016) described counting and measurement methods practised by the Kabihug people of sitio Calibigaho, Jose Panganiban Camarines Norte, Philippines. Rosin (1984) described how a non-literate peasant, Rupsingh, demonstrated his capabilities in arithmetical computation by counting finger joints and in estimating the population of cities such as Jaipur and Calcutta to high accuracy. Botetano and Abrahamson (2022) make use of body-based arithmetic methods to engage Indigenous students to learn addition and multiplication.

*Architecture*

Geometric and mensuration methods used in Indigenous architecture are the subjects of 19 of the studies (Supplement 3). Omthuan et al. (2009) derived numerical relationships used in Isan house building. The non-standard, or 'inaccurate', method of taking measurements is a purposeful way of creating a bespoke relationship between the house and the owner. Measurements are not based on standard units (e.g., inches) but on the owner's body components such as thumb, finger, arm and hand (or anthropometry). Darmayasa et al. (2018) found that the concept of linear regression has been applied in constructing the *Bale Saka Roras*, a Balinese traditional house. The formula

$$Y = a + bX + cZ$$

is used to estimate the length of pillars, where $Y$ is the dependent variable representing the length of the pillar, $a$, $b$, $c$ are constants, $X$ represents the length of house owner's index finger, and $Z$ represents the length of house owner's middle finger.

Supiyati and Hanum (2019) examined the *sasak* granary building technology in Lombok, which dates from 3500 years BCE, where the size of buildings is determined similarly by anthropometry, and identified possible relations between the design of the roof of a *sasak* building and a sinusoidal function. Hermanto and Nurlaelah (2019) describe how birth dates are used by the Kampung Naga people to determine the construction area of the house built for a newly married couple, based on the formula

$$A = Vh \times Vw,$$

where $A$ is the area of the house to be built, $Vh$ is a length value determined by the day of the husband's birth and $Vw$ is a length value determined by the day of the wife's birth. Each day of the week in the Kampung Naga calendar is linked to a specific number, e.g., Monday to 4, Tuesday to 3, and Wednesday to 7. The area $A$ acquired from the calculation is divided by 3.

If the remainder is 0, the final construction area needs to be subtracted by 1 and the lengths adjusted as necessary.

*Ceremonies and divinations*

Activities in cultural ceremonies often require skills in arithmetic, probability, and pattern recognition (9.0% of the survey literature, in Supplement 3). Budi Darmayasa et al. (2018) identified that integer, least common multiple, and mixed fraction operations are used to choose festival days in the Balinese Hindu calendar system.

Divination in various forms has long been a popular activity in human societies (Peek, 1991), and is still to this day. In 1937, Evans-Pritchard observed the divination practices of the Azande people in Central Africa and noted that diviners consider multiple possible outcomes of a divination and assign different levels of likelihood to each one (Eglash, 1997).

Application of the modulo operation is ubiquitous in divination systems throughout Africa, such as the *ifa* of Yoruba communities (Abimbola, 2008; Morton-Williams et al., 1966). Ascher (1997) studied *sikidy*, a method of divination used in Madagascar, and recognized that the mathematics of *sikidy* is analogous to symbolic logic and makes use of parity checking. Chemillier (2009) undertook a more detailed, multi-disciplinary collaborative investigation of this Malagasy divination system. He observed that some master diviners could work the system in a predictive manner, and likened their understanding and mental skills to those of a chess grandmaster visualizing several moves ahead.

Indonesian researchers on ethno-mathematics have also studied traditional divination systems. Kusuma et al. (2017) described how the Cipatujah people in West Java use subtraction, division, and remainders to calculate *repok*, a numerical result showing the compatibility of a man and woman, to determine whether they could develop harmonic relationships in marriage. Budi Darmayasa et al. (2018) revealed how Bali Mula people in Kintamani District use least common multiples and mixed fractions to make their *Wuku* calendar and schedule festivals. Prahmana et al. (2021) translated the process of selecting dates for death ceremonies in Yogyakarta into the formula $H = c + b \pmod 7 - 1$, where $H$ is the date of the funeral ceremony, $c$ is the number of the day of death and $b$ is the mourn day, the value of which can be 3, 7, 40, 100 and 1,000 according to the local culture. U. Umbara et al. (2021) mathematized the approach that Indigenous communities in Cigugur use to determine the best day to start constructing a house. Utami and Sayuti (2019) analyzed the Javanese tradition of calculating *weton*, a quantity used in traditional matchmaking and marriage customs in Java. The success or harmony of a proposed marriage is predicted by comparing the score or ranking of the couple's *weton*. Utami and Sayuti (2019) transformed the procedure to numerical remainder formulas.

### *Games and puzzles*

Indigenous or traditional societies master probability, statistics and combinatorics through games and puzzles (10.3% of the survey literature, in Supplement 3). Stillman (1995) documented the combinatorics of the game 'Pitch and Toss' played by the Igbo people of Nigeria. Turmudi et al. (2021) described the arithmetic operations used in the Malang traditional game Tong Tong Galitong Ji, in which players are required to carry out addition, multiplication and subtraction accurately using their fingers. Bjarnadóttir (2009) reveals the ancient roots of arithmetic and algebraic reasoning of the thirty-bird puzzle, which was known in ancient Chinese and Indian societies. This puzzle is an early exemplar of an underdetermined system, for which in general there are infinitely many solutions, or none. A problem recorded in the treatise from c. 485 CE known as The Mathematical Classic of Zhang Qiujian treats such a system: 'Now one cock is worth 5 qian, one hen 3 qian and 3 chicks 1 qian. It is required to buy 100 fowls with 100 qian. In each case, find the number of cocks, hens and chicks bought'(Bjarnadóttir, 2009). Possible solutions recorded in the book are: 4 cocks, 18 hens and 78 chicks; 8 cocks, 11 hens and 81 chicks; 12 cock, 4 hens and 84 chicks. In algebraic form the problem is written as $x + y + z = 100$, $5x + 3y + z/3 = 100$. Tiennot (2023) documents a sowing game, known as solo, in Tanzania, Comoros, Madagascar and Mozambique, and constructs ethno-models of combinatorial strategies of experts in playing the game.

**Semiotics**

Mathematical language is supposed to be a unified semiotic system, adjunct to natural language, for expressing its worlds of objects meaningfully. Since mathematics is entirely a human construction, the 'truths' it establishes are attributes of the 'mathematicians' who read and write mathematical symbols and accept its persuasions (Rotman, 1988). European mathematics developed its current semiotic system internally and from Greek, Hindu, Arabic, and other external sources, but there are semiotic systems in non-European cultures.

Vergani (1999) studied the symbolism and meanings of the numbers used by the Dogon, an ethnic group in Mali. The Dogon have two language dialects, one of which (*Sigi So*) is a secret, rhythmic language different from their daily tongue. *Sigi So* mythology seems to describe nonlinear growth, and preliminary studies suggest the foundational narratives may follow a power law. Manapat (2011) found that pre-colonial Tagala people in Philippines had well-developed enumeration and arithmetical systems before colonization by the Spanish. The system of counting in the Tagalog language names numbers as large as a hundred million. Muhtadi (2019) elucidated the mathematical operations in Sundanese manuscripts from West Java, including addition, subtraction, multiplication, division, modulo and set theory. Yulianto et al. (2020) decoded the mathematics in Dhikr, a form of Islamic meditation having different ritualistic symbols, notably the number 165. Practitioners can choose the number of times they pray, $N$, subject to two conditions: $N \geqslant 165$ and $N = 2k + 1$, where $k$ is any natural number. The Arabic letters of the phrase 'laa ilaaha illallah' each is represented by its corresponding number from the Abjad numerals system, which gives 12 distinct numerical values. Summing all these values gives the total of 165.

*Handicraft*

Non-literate handicraft artisans are included in five studies as a group that develops and translates mathematical thought into practical skills. Millroy (1991) described the geometric ideas, including the concepts of straightness, parallel lines, axes of symmetry and right angles, used by a group of nonliterate carpenters in Cape Town, South Africa. Masingila (1994) observed a group of carpet installers and described the trade-specific, non-school-based mathematics of measurement, computational algorithms, geometry and proportion applied by the installers. Novitasari et al. (2023) described the use of two- and three-dimensional shapes and principles of transformation and dilation in creating the Sasaknese traditional instrument, *gendang beleq*. Prahmana and Istiandaru (2021) identified set theory in Javanese shadow puppetry, by the categorization of puppet characters into distinct sets, such as 'Simpingan kiri' for those with malevolent traits and 'Simpingan kanan' for benevolent ones. During the show, characters may transition between sets, e.g., a traditionally 'bad' character assuming a guise of virtue, or vice versa.

*Navigation*

Indigenous navigation knowledge includes advanced spatial reasoning, understanding of celestial periodicities, and pattern recognition from natural cues. Mathematical abstractions and modelling were identified by Ascher (1995) in her study of the stick charts, or *rebbelith*, made by the Marshall Islanders to aid navigation of the Pacific Ocean, where there are long and fast moving swells, and for the purpose of teaching navigation.

Abdullah (2017) created an ethnomodel of the tide prediction method used by the Sundanese people of West Java (see flow chart of Figure S4 in Supplement 1). The timings of the tides are calculated by referring to the *hijriyah* lunar calendar. The Sundanese calibrated the arrival time of the first tide from 12 am each day at Santolo Beach to correspond with the *hijriyah* dates. Arrival times of the 2 maxima and 2 minima could be calculated from the two formulas in the model: The time of the first high tide in the next 24 hours $= \frac{4}{5} \times$ the date in the hijriyah calendar; and the hour of the second tide = the hour of the first tide + 6 hours.

*Music*

Four of the studies present the rhythms, patterns, symmetries, and sequences of traditional and Indigenous dance and music. Indrawati et al. (2021) analyzed the fractional patterns in notated Surabaya regional songs. Radiusman et al. (2021) described intricate geometry in the hand gestures of Bali's Pendet Dance, in which dancers employ acute, right, and obtuse angles in named movements such as *ngumbang*, *agem*, *ulap-ulap*, and *ngelung*.

*Agriculture*

Four of the studies document traditional or Indigenous agricultural practices that incorporate mathematical principles in crop spacing, irrigation calculations, and harvest predictions. Suprayo et al. (2019) recorded local area measurement and mass calibration methods used by Suranenggala Kidul villagers in Cirebon.

*Trade*

Two of the studies elucidate mathematical strategies used for bartering, estimation, and value comparison. Nurjanah et al. (2021) documented a calculation system used in the Marosok trading tradition by the Minangkabau tribe in West Sumatra, by which addition and subtraction of 3-digit numbers was carried out using only the hands.

*Kinship*

Indigenous kinship systems can be analyzed and understood in terms of group theory. Ascher and Ascher (1986) display the connections between group theory and kinship relationships by studying the regulation of marriage among the Aranda people of Australia and the Armbrym people of Vanuatu.

*Land surveying*

Traditional or Indigenous land surveying techniques show accurate spatial understanding of terrains and boundaries. Williams and Jorge y Jorge (2008) and Jorge et al. (2011) assessed the accuracy of Aztec land surveys recorded in the Codex Vergara, also known as 'cadastral' documents, that contain detailed measurements of plots of land. Errors of the Acolhua-Aztecs land survey records made in the 1540s were found to be of similar magnitude to those made by English surveyors a century later and even lower than those of American surveyors two or three centuries later.

*Calendars*

Three of the studies document arithmetical sequences and modular arithmetic used in calendars and chronograms to chart cyclical patterns of lunar and solar movements to high accuracy. Syahrin et al. (2016) found that divisions and remainders are broadly used in the Indigenous Aboge calendar of West Java to determine festivals and important dates.

*Artifacts*

Three of the studies are on traditional and Indigenous artifacts that manifest mathematical concepts such as geometric patterns and symmetries and/or have a mathematical purpose, such as algorithm and numeral system 'hard copy'. Ascher and Ascher (2013) identified a *quipu* that is likely the earliest record of the 'grating method'. This multiplication technique was previously believed to have been developed in 15th C Europe.

*Sand drawing*

In two of the studies the seemingly arcane activity of sand drawing is the subject of attention. Early work from M. Ascher (1988) introduces the graph theory and practice of the Malekulan cultures of Vanuatu, where the art of drawing continuous line figures in sand, or *nitus*, is an essential element of important rituals. Each single-line figure should be finger-drawn in continuous motion, such that the finger is not lifted from the sand until the figure is completed at the starting point. Figure S5 in Supplement 1 shows four variations of a recorded basic *nitus*, used by the Indigenous people of Malekula, that transmit qualitatively similar information as they have the same graphical properties.

The sand drawings (*sona*) of the Tchokwe people of the west central Bantu area of Africa are integral to storytelling and education. They can be analyzed in terms of the greatest common divisor (GCD) of integers (Demaine et al., 2007), and present geometric challenges such as finding the minimum-length *sona* map for a set of dots (Barros & da Silva, 2022). Demaine et al. (2007) further elucidated the connection between the Tshokwe *sona* and graphs supposedly invented by Gauss in 1839. They proposed a new approach to solve the famous Traveling Salesman Problem inspired by the minimum-length *sona* figures.

*Food preparation*

Traditional or Indigenous food preparation can embody mathematical principles through precise proportions, symmetrical arrangements, and the calculation of cooking durations and sequences. Pathuddin and Ichsan Nawawi (2021) found that the process of making *Barongko* cakes by the Buginese community uses division, congruence and similarity, as well as 3-dimensional shapes such as triangular prisms and half spheres.

*Fishing*

Fish-catching methods depend on physical and mathematical understanding of water flows and landscape geometry. Malalina et al. (2020) documented fishing activities in the Musi River that use Cartesian coordinates for optimal location selection, advanced measurements of weight and velocity, and geometric concepts in the design of fishing tools.

## 4. Discussion

Babbitt et al. (2012) found that ethnomathematics and Indigenous mathematics research face two challenges. First, researchers must investigate mathematical ideas in cultural practices that are often assumed to be unrelated to mathematics, with 'mathematics' by default taken to mean 'European mathematics'. A question frequently asked by educators seeking to include Indigenous content and perspectives in mathematical education is 'Where should I look for Indigenous mathematics?' Second, even if researchers can identify previously unrecognized mathematics, it is difficult to decode and present it in a 'formal' approach endorsed by mainstream mathematicians and apply this knowledge.

This scoping review unearths a third, more concealed or perhaps deliberately overlooked challenge: the gatekeeping problem. While the first two challenges focus on the identification and translation of mathematical concepts within cultural contexts, this third challenge refers to the more insidious institutional, social and cultural barriers. Strongly held and propagated views in sections of the mainstream mathematics community deny the existence of Indigenous mathematics capabilities, or, if confronted with evidence, claim European priority or that mathematical knowledge was taught by outsiders from Europe (Deakin, 2010).

The gatekeeping challenge is linked to another problem: the dearth of Indigenous mathematics researchers. This has manifested in a discernible 'white saviour' undertone within the field, where Western scholars tend to dominate the exploration of Indigenous knowledge systems, often without ensuring equitable collaborations or, in some instances, bypassing partnerships altogether. Progression and exploration of Indigenous mathematics must be led by Indigenous voices. However, this path is riddled with inherent contradictions. While it is essential for Indigenous scholars to reclaim their mathematical traditions, they are simultaneously expected to navigate and master a curriculum deeply rooted in European-based mathematics. This juxtaposition underscores the need for a balanced approach and emphasizes the importance of addressing the gatekeeping issues to ensure the rightful representation and leadership of Indigenous scholars in the domain of ethnomathematics.

Those challenges also act as barriers to constructive communication between ethno-mathematicians and mainstream mathematicians. Although cultural mathematical knowledge and practices have been well-documented (de Almeida, 2019; Goetzfridt, 2007), pushback against the evidence also persists. Ethnomathematics has been viewed as a hindrance to the progress of mainstream mathematics as pure, context-free and universal truths (Ernest, 1994). An unfortunate unintended effect has been to create a widening rift between ethno-mathematicians and mainstream mathematicians (Skovsmose, 2022).

Ethnomathematics promotes the understanding that mathematical knowledge is not the exclusive domain of Western academia (Ascher & Ascher, 1986). Its research outputs reveal that diverse cultures around the world have developed sophisticated mathematical systems tailored to their particular contexts and needs. By studying these systems, we can gain insights into various ways of knowing and interpreting the world mathematically (Ascher & Ascher, 1986). In contrast, European mathematics, with its rapid and specialized advancements, tends to greater detachment from worldly applications (Kline, 1982). While this view is debated among scholars, as discussed in Section 3.2, Indigenous or traditional mathematical systems usually are tied to real-world contexts – from navigation to agriculture, from ceremony to trade. Such systems have evolved out of immediate necessity and thus carry with them historical and contemporary connections to everyday life. Mainstream mathematicians tend to view the culture-centric work of ethno-mathematicians as closer to anthropology, history, archeology or philosophy (Barton, 2007), having little relevance to 'real mathematics'. This estrangement – which we hope is temporary – may at least give space and time for both groups of mathematicians to grow up. Ethnomathematics itself has developed different streams of research. Based on bibliographic mapping of the co-occurrence of keywords used by the authors and bibliographic indexes, we have found two clusters of research. These are visualized in the graphic of Figure 3.

The right cluster contains research before 2016 including most of the studies in Stages 1 and 2 of Figure 2, connecting to the most recent research on ethno-modelling (co-occurrences

= 4; total link strength = 12) and glocalization[3] (co-occurrences = 2; total link strength = 8). Some studies in this right-hand cluster focus on bringing mainstream mathematical concepts, such as graph theory, group theory and statistics, into an ethnomathematical knowledge system. Others explore how Indigenous mathematics may be decoded from cultural practices and become more 'scientifically' accepted, and the potential benefits of studying Indigenous or traditional mathematics.

---

[3] The term 'glocalization', which was coined by Orey, D. C., & Rosa, M. (2021). Ethnomodelling as a glocalization process of mathematical practices through cultural dynamism [Article]. *Mathematics Enthusiast*, *18*(3), 439-468. https://www.scopus.com/inward/record.uri?eid=2-s2.0-85101544170&partnerID=40&md5=a1dc2adcaeb35803fb0d28dc654ea7e1 , reflects the need for integration and acknowledgement of local and Indigenous mathematical knowledge within a global understanding of mathematics.

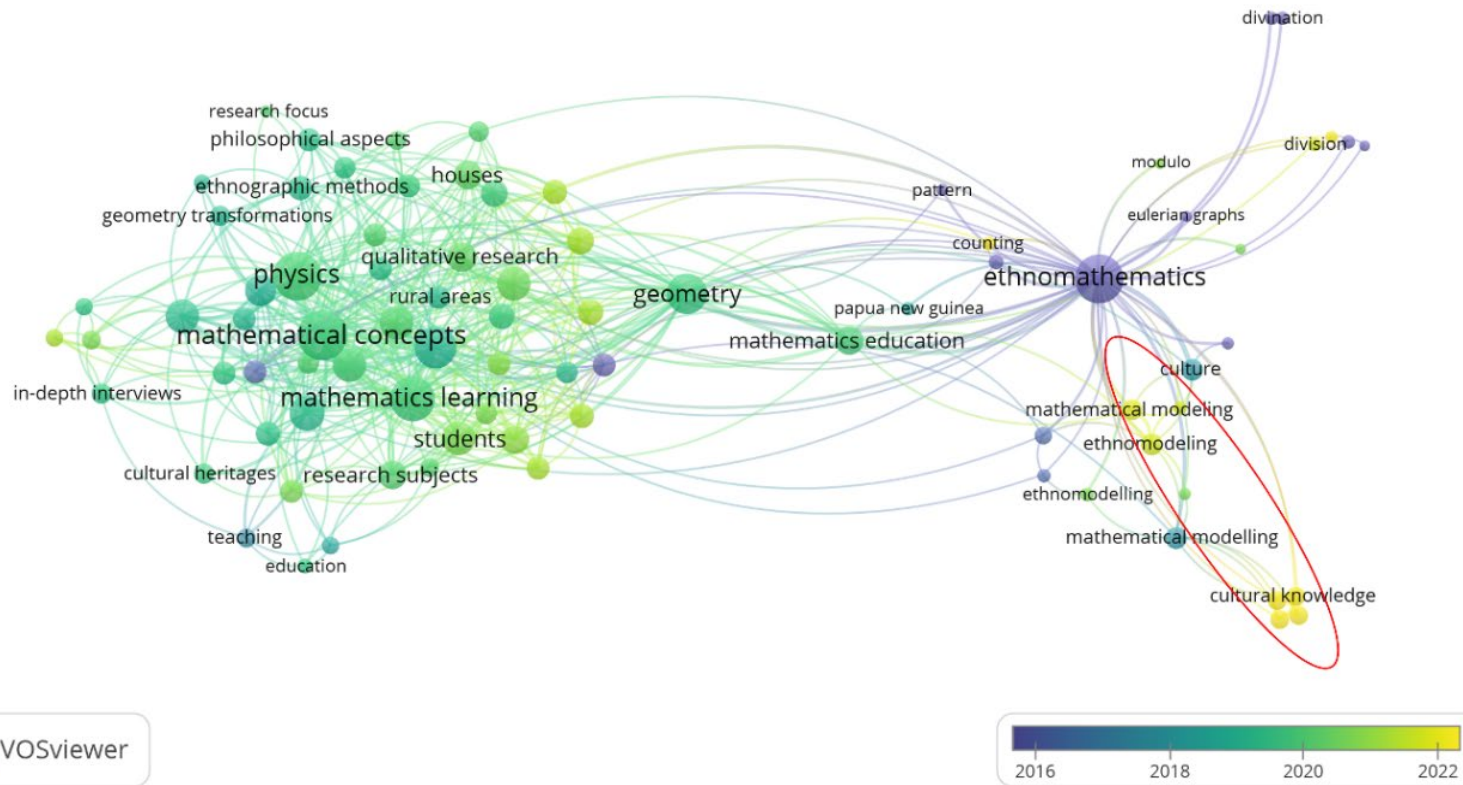

Figure 3. Network analysis diagram, where the nodes (or vertices) represent keywords, with type size proportional to their relative frequencies, and the or edges link co-occurrences of keywords within the same study. The thickness of a link (total link strength) indicates the frequency of such co-occurrences. A keyword connected by multiple thick links suggests its central importance in the field's themes. The most recent research is in the red ellipse. For precise co-occurrences between keywords and link strengths of each keyword, refer to the Supplement 4. Others explore how Indigenous mathematics may be decoded from cultural practices and become more 'scientifically' accepted, and the potential benefits of studying Indigenous or traditional mathematics.

The left cluster of Figure 3 consists of studies that concentrate on elucidating mathematical knowledges (co-occurrences = 19; total link strength = 87) of Indigenous and traditional communities, and investigating how that could be used to improve mathematics learning (keyword occurrences = 13; total link strength = 66). Most of the recent research in this cluster explores mathematical ideas and concepts across different cultural groups. These studies have received minimal attention from mainstream mathematics after 2006, but this research avenue empowers Indigenous knowledge holders to inherit and re-visit their knowledges and opens a pathway for outsiders to gain a deeper appreciation and understanding of the mathematical concepts and philosophies embedded in Indigenous cultures.

Although the two challenges proposed by Babbitt et al. (2012) still exist, the booming studies of cultural mathematics from 2012 evince the growth and vitality of the field. In Figure 3, research results from the left cluster provide guidance on where to find Indigenous mathematics, and those in the right cluster are complementary in decoding methodologies. Results from both clusters of studies are needed to progress research on Indigenous mathematics, empower Indigenous knowledge holders and decolonize mathematics. To extend the influence of ethnomathematics, more research is needed on decoding and systematizing mathematical knowledges identified, without damaging the cultural heritage.

Research on mathematics of marginalized groups will continue to yield results that may surprise the gatekeepers of European-based mathematics. As shown in this scoping study, different cultural groups, historically and contemporaneously, have chosen or choose to use and master mathematics based on their specific needs. This does not mean that mathematics in traditional and Indigenous communities cannot be broadly applied (Thomas, 1996). Research into Indigenous mathematics and ethnomathematics serves as a bridge that connects knowledge systems, with potential to offer an enriching palette of new methodologies and epistemologies. When we cast a wider net, we recognize mathematics, far from being monolithic, has evolved in many forms. In contrast to the esoteric abstractions of contemporary mainstream mathematics (Kline, 1982), Indigenous or traditional systems frequently exhibit direct applicability, from sustainable agriculture to clever architectural designs. This relevance has considerable pedagogical potential, offering culturally responsive teaching that can rekindle interest in mathematics among communities that have felt historically alienated from mainstream curricula (Meaney et al., 2021).

A comprehensive exploration of non-European mathematics may open doors to holistic solutions for contemporary challenges. Cultural burning, as practised by Indigenous Australian communities, is a part of a holistic approach to land management that incorporates fire as a tool (Neale, 2022). The principles of cultural burning are rooted in mathematical understanding of biological, ecological, and planetary cycles and its application uses understanding of reactive fluid flows. This holistic approach surpasses the scope of mainstream disciplinary science and mathematics. Indigenous knowledge keepers can identify the right time to burn based on a multitude of signals from the environment, such as the flowering of specific plants or the behaviors of certain animals. In recent years, as the intensity and frequency of bushfires have increased, we have seen wide acceptance of the importance of Indigenous fire management practices.

While ethnomathematics is concerned with the mathematical practices of diverse cultures, Indigenous mathematics is more specific. The former champions the notion that mathematical thinking is universal but manifests differently across cultures, while the latter stresses the unique mathematical systems rooted in Indigenous traditions. This dichotomy is evident in the bibliographic mapping in Figure 3. The right cluster, mostly pre-dating 2016, leans towards synthesizing Indigenous mathematics with more widely accepted mathematical systems, aiming for a broader academic acceptance. In the left cluster we see more identification and appreciation of mathematical concepts as they manifest in various cultural contexts. While there is undeniable value in understanding Indigenous mathematical concepts on their own terms, we also need collaborative efforts where Indigenous scholars play leading roles, ensuring the authenticity of representation, and which Western-trained mathematicians approach with open-mindedness and respect, ready to learn and expand their horizons. The aim would be to create a pluralistic but holistic mathematical landscape, where various systems coexist, collaborate, and enrich one another.

**Conclusion**

In Western societies mathematics is seen as an elite branch of knowledge, the provenance of which is almost exclusively European. Such a perspective has perpetuated the myth that mathematical proficiency is poor or absent in traditional and Indigenous communities. However, as elucidated in this scoping study, these communities have always harbored rich mathematical knowledge that is intertwined with their cultural practices and worldviews, challenging conventional notions of what constitutes 'mathematics'. To acknowledge and embrace these diverse mathematical expressions is not merely an act of respect towards these communities but may open new opportunities for developing new mathematics and sciences. Innovations such as ethno-modelling and ethno-computing are tapping into this potential. As the field of ethnomathematics burgeons, it is imperative to navigate it with intention and depth. It is not enough simply to scout for mathematical elements across diverse cultures. A profound, respectful engagement is crucial not just for unearthing the latent mathematical wisdom these communities hold, but also for quenching gatekeeper critiques that deny the existence and validity of Indigenous mathematics. Such sincere engagement directly challenges these dismissive stances, by emphasizing that mathematics across cultures is not differentiated by the knowledge itself but by its presentation, application, and transmission. By recognizing this, we validate diverse knowledge systems and invite an enriched, collaborative evolution of the mathematical sciences. The journey ahead in ethnomathematics is promising in terms of new mathematical knowledge that potentially may seed innovations that link theory to tangible outcomes.

**Supplementary Document 1**

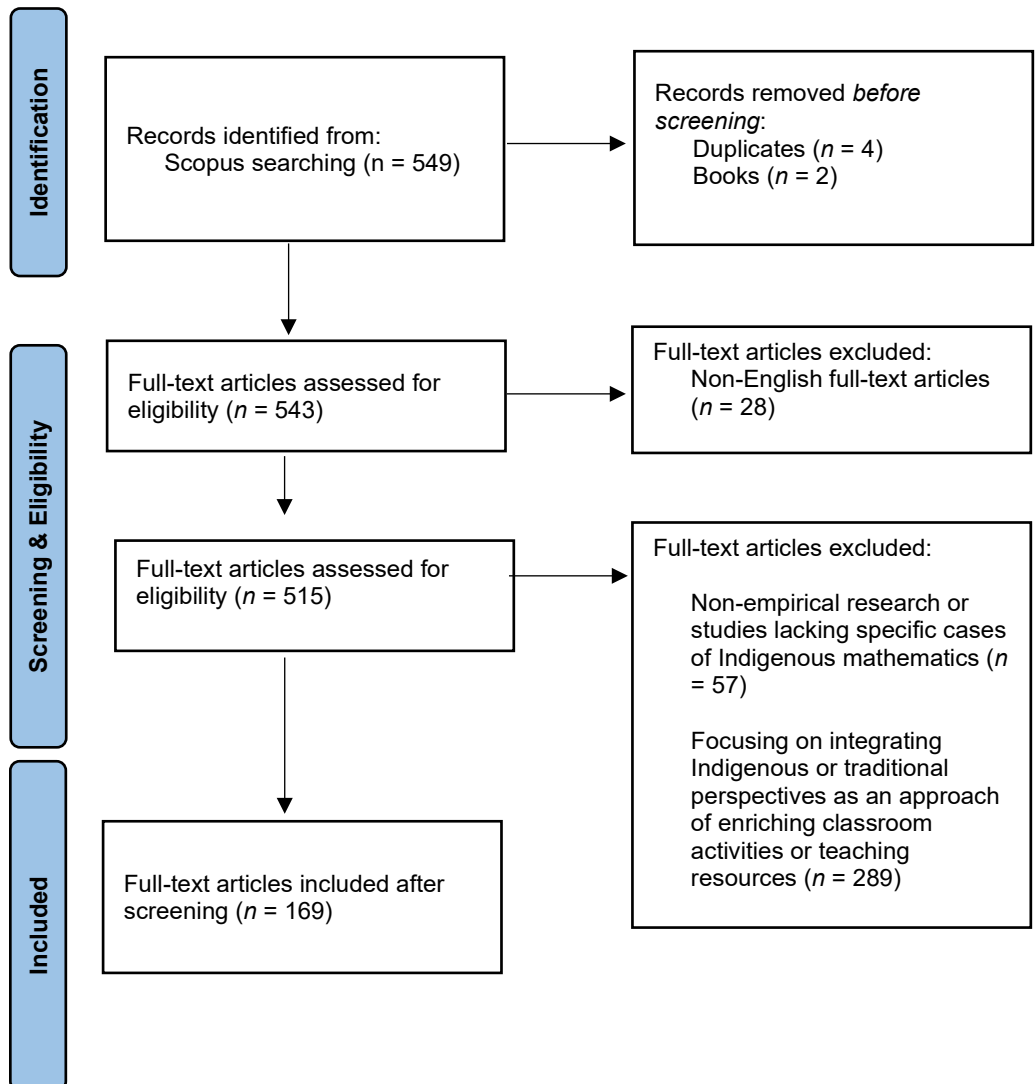

Figure S1. Flow diagram of the literature search strategy and review process described in Section 2 of the main article

Figure S2. Word cloud displaying 200 most frequent words in the abstracts of the 169 included studies (produced by Nvivo 12). The 5 most frequent words are: mathematics (559) > culture (277) > ethnomathematics (215) > concepts (184) > tradition (108)

Table S1. Example of an algebraic weaving pattern based on $(a + b + c + d + e + f)^2$. The pattern below starts from left to the right and each column specifies an interaction between two threads. The 36 terms of the expansion are enumerated across the top and (the first 6) down the left side of the draft's music-like notation bars.

|   | 1 | 2 | 3 | 4 | 5 | 6 | 7 | 8 | 9 | 10 | 11 | 12 | 13 | 14 | 15 | 16 | 17 | 18 | 19 | 20 | 21 | 22 | 23 | 24 | 25 | 26 | 27 | 28 | 29 | 30 | 31 | 32 | 33 | 34 | 35 | 36 |
|---|---|---|---|---|---|---|---|---|---|---|---|---|---|---|---|---|---|---|---|---|---|---|---|---|---|---|---|---|---|---|---|---|---|---|---|---|
| 1 | $a$ | $a$ | $a$ | $a$ | $a$ | $a$ | $a$ | $a$ | $a$ | $a$ | | | | | | | | | | | | | | | | | | | | | | | | | | |
| 2 | | $b$ | $b$ | | | | | | | | $b$ | $b$ | $b$ | $b$ | $b$ | $b$ | $b$ | $b$ | $b$ | | | | | | | | | | | | | | | | | |
| 3 | | | $c$ | $c$ | | | | | | | | $c$ | $c$ | | | | | | | $c$ | $c$ | $c$ | $c$ | $c$ | $c$ | $c$ | | | | | | | | | | |
| 4 | | | | $d$ | $d$ | | | | | | | | $d$ | $d$ | | | | | | | $d$ | $d$ | | | | | $d$ | $d$ | $d$ | $d$ | $d$ | | | | | |
| 5 | | | | | $e$ | $e$ | | | | | | | | $e$ | $e$ | | | | | | | $e$ | $e$ | | | | | $e$ | $e$ | | | $e$ | $e$ | $e$ | | |
| 6 | | | | | | $f$ | $f$ | | | | | | | | $f$ | $f$ | | | | | | | $f$ | $f$ | | | | | $f$ | $f$ | | | | $f$ | $f$ | $f$ |

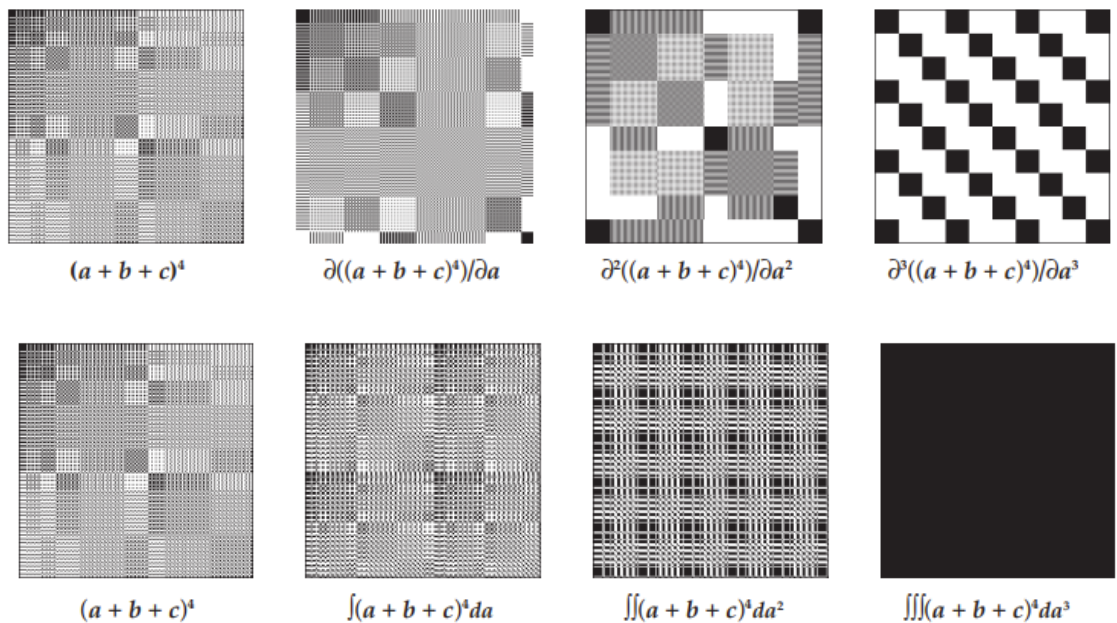

Figure S3. Renderings of computed weaving patterns generated by the polynomial expression $(a + b + c)^4$, along with its first, second and third partial derivative with respect to $a$ (upper); and its first, second and third integral with respect to $a$ (lower), where integration constants are not given but presumably are set to zero. (Reproduced from Griswold, 2002, pp. 9–10)

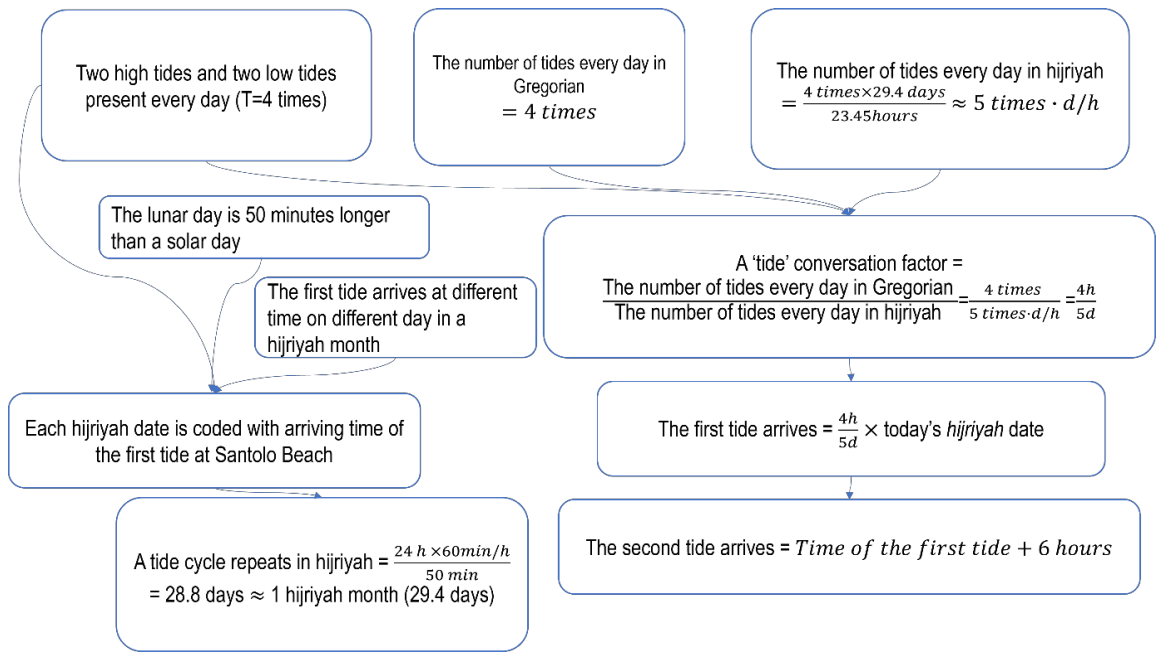

Figure S4. An ethnomathematical model of Sundanese knowledge on local tide cycles and prediction. (Abdullah, 2017).

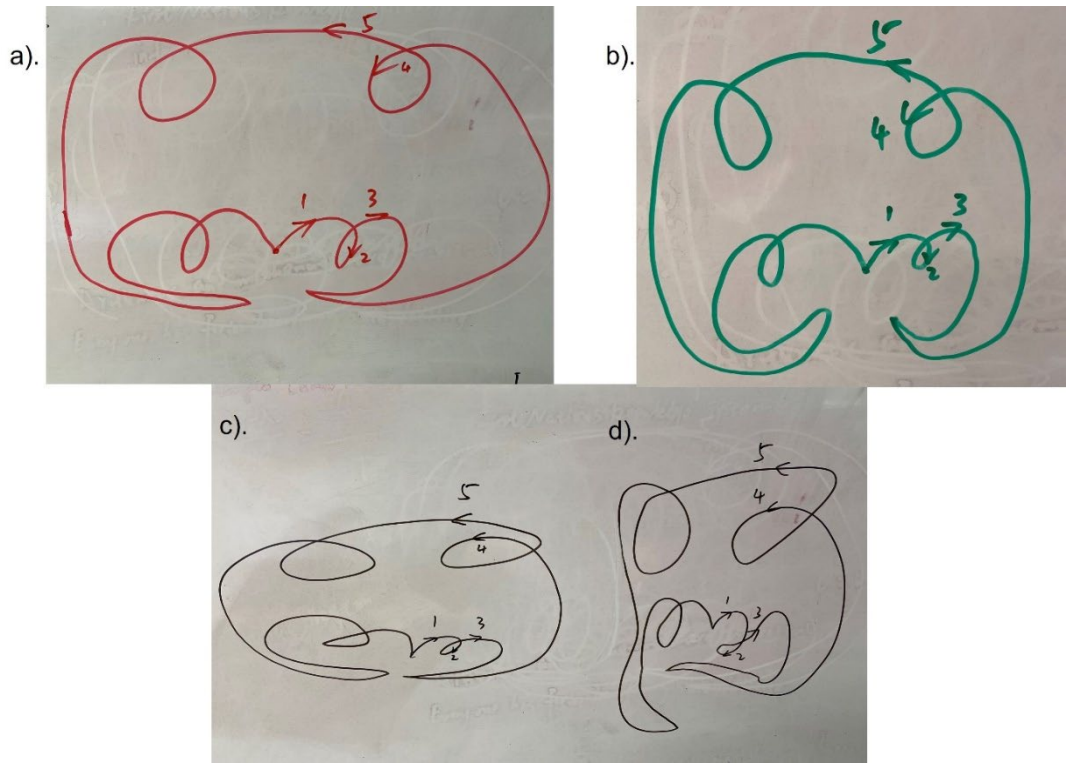

Figure S5. Four variations of a *nitus* sharing the same graphical features